\documentclass[12pt]{amsart}

\newcommand {\supplus}{\mathop{{\supset}\llap{\raise 
0.5pt\hbox{\normalfont\small+}\hskip 0.5pt}}} 

\newcommand {\subplus}{\mathop{{\subset}\llap{\raise 
0.5pt\hbox{\normalfont\small+}\hskip 0.5pt}}}  

\newcommand {\Cee}    {{\mathbb  C}}

\newcommand {\Zee}    {{\mathbb  Z}}

\newcommand {\fA}     {{\mathfrak{A}}}

\newcommand {\fg}     {{\mathfrak{g}}}    %
\newcommand {\fgl}    {{\mathfrak{gl}}}  %
\newcommand {\fh}     {{\mathfrak{h}}}

\newcommand {\fm}     {{\mathfrak{m}}}

\newcommand {\fosp}   {{\mathfrak{osp}}}

\newcommand {\fsl}    {{\mathfrak{sl}}}

\def \opname#1#2%
  {\expandafter\newcommand \csname #1\endcsname {{\mathop{#2}\nolimits}}}

\newcommand{\rmname}[1]
  {\expandafter\newcommand \csname #1\endcsname {{\operatorname{#1}}}}

\newcommand{\rmnameii}[2]
  {\expandafter\newcommand \csname #1\endcsname {{\operatorname{#2}}}}

\rmname{act}
\rmname{Ad}
\rmname{Add}
\rmname{ad}
\rmname{Alt}
\rmname{alt}
\rmname{Ann}
\rmname{antidiag}
\rmname{Ber}
\rmname{ber}
\rmname{Br}
\rmname{card}
\rmname{ch}
\rmname{Char}
\rmname{cem}
\rmname{cj}
\rmname{Cliff}
\rmname{cntr}
\rmname{codim}
\rmname{coind}
\rmname{const}
\rmname{col}
\rmname{cork}
\rmname{cpr}
\rmname{diag}
\rmnameii{Div}{div}
\rmname{Def}
\rmname{Der}
\rmname{Dim}
\rmname{End}
\rmname{Even}
\rmname{Ext}
\rmname{gr}
\rmname{Hom}
\rmname{HT}
\rmnameii{Ht}{ht}
\rmname{hwt}
\rmname{Id}
\opname{IM}{Im}
\rmname{id}
\rmname{ind}
\rmname{Ind}
\rmname{Inf}
\rmname{irr}
\rmname{Le}
\rmname{Lie}
\rmname{lwt}
\rmname{mult}
\rmname{Mor}
\rmname{nm}
\rmname{Ob}
\rmname{Odd}
\rmname{Osc}
\rmname{per}
\rmname{Pic}
\rmname{pr}
\rmname{pro}
\rmname{Prime}
\rmname{Proj}
\rmname{prt}
\rmname{pt}
\rmname{Q}
\rmname{qet}
\rmname{qtr}
\rmname{rd}
\rmname{rk}
\rmname{row}
\rmname{Res}
\rmname{salt}
\rmname{Sch}
\rmname{SBr}
\rmname{scalar}
\rmname{Ser}
\rmname{sign}
\rmname{Smbl}
\rmname{spin}
\rmname{ssym}
\rmname{str}
\rmname{st}
\rmname{sgn}
\rmname{sq}
\rmname{symm}
\rmname{supp}
\rmname{Supp}
\rmname{St}
\rmname{Spec}
\rmname{Spm}
\rmname{tr}
\rmname{vpt}
\rmname{weyl}
\rmname{Weyl}
\rmname{Witt}

\opname{vvol}  {{v\hspace{-0.1ex}o\hspace{-0.02ex}l\/}}
\opname{pnt}  {\text{\normalfont pt}}
\opname{Span} {{Span}}
\opname{slim} {\overline{\lim}}
\opname{Vol}  {{V\hspace{-0.55ex}o\hspace{-0.02ex}l\/}}
\opname{Par}  {{P\hspace{-0.3ex}a\hspace{-0.05ex}r\/}}

\rmname{Mat}
\rmname{Bil}
\rmname{Diff}
\rmname{Ker}
\rmname{Herm}
\rmname{Coker}
\rmname{Conn}
\rmname{Covect}
\rmname{Vect}
\rmname{Int}

\newcommand {\eps} {\varepsilon}

\newcommand {\tto} {\longrightarrow}
\newcommand {\pder}[1] {{\frac{\partial}{\partial {#1}}}}
\newcommand {\pderf}[2] {{\frac{\partial {#1}}{\partial {#2}}}}

\newcommand {\bcdot}   {\mathbin{\hbox{\raise.4ex\hbox{\bf.}}}} 

\newcommand {\secno} {}
\newcommand {\ssecfont} {\normalfont\bf}

\newtheorem{Theorem}{\secno Theorem}

\newtheorem{Lemma}[Theorem]{\secno Lemma}

\newenvironment {th*}[1]
    {\gdef\thname{#1} \begin{thn}}%
    {\end{thn}}
\newtheorem{thn}[Theorem] {\thname}

\theoremstyle{definition}

\newenvironment {ex*}[1]
    {\gdef\thname{#1} \begin{exn}}%
    {\end{exn}}
\newtheorem{exn}[Theorem]{\thname}

\theoremstyle{remark}

\newtheorem{Remark}[Theorem]{\secno Remark}

\newenvironment {rem*}[1]
    {\gdef\thname{#1} \begin{remn}}%
    {\end{remn}}
\newtheorem{remn}[Theorem]{\thname}

\newcommand {\ssbegin}[2]
  {\def \secno {\gdef \secno {}{\ssecfont #1. }}%
   \begin{#2}}
\newcommand {\ssec}{\subsection*}

\begin{document}

\title[Enveloping algebra and orthogonal
polynomials]{Enveloping algebra $U(\fgl(3))$ and orthogonal
polynomials in several discrete indeterminates}

\author{A.~Sergeev}

\address{(Correspondence) Department of Mathematics, University of
Stockholm, Kr\"aftriket hus 6, S-106 91, Stockholm,
Sweden; mleites@matematik.su.se (On leave of absence from Balakovo
Institute of Technique of Technology and Control, Branch of Saratov
Technical University, Balakovo, Saratov Region, Russia)}

\thanks{I am thankful to D.~Leites for encouragement and help and to 
ESI, Vienna, for hospitality and support.}



\begin{abstract} Let $A$ be an associative algebra over $\Cee$ and $L$ 
an invariant linear functional on it (trace).  Let $\omega$ be an 
involutive
antiautomorphism of $A$ such that $L(\omega(a))=L(a)$ for any $a\in A$.  
Then $A$ admits a symmetric invariant bilinear form 
$\langle a, b\rangle=L(a\omega(b))$.  For $A=U(\fsl(2))/\fm$, where 
$\fm$ is any maximal ideal of $U(\fsl(2))$, Leites and I have 
constructed orthogonal basis whose elements turned out to be, 
essentially, Chebyshev and Hahn polynomials in one discrete variable.

Here I take $A=U(\fgl(3))/\fm$ for the maximal ideals $\fm$ which 
annihilate irreducible highest weight $\fgl(3)$-modules of particular 
form (generalizations of symmetric powers of the identity 
representation). In whis way we obtain multivariable analogs of Hahn 
polynomials.
\end{abstract}

\maketitle

This paper appeared in: Duplij~S., Wess~J. (eds.)  {\em Noncommutative
structures in mathematics and physics}, Proc.  NATO Advanced Reserch
Workshop, Kiev, 2000.  Kluwer, 2001, 113--124; I just want to make it
more accessible.

\section*{\S 1. Background}

\ssbegin{1.1}{Lemma} Let $A$ be an associative algebra generated 
by a set $X$.  Denote by $[X, A]$ the set of linear combinations of 
the form $\sum [x_{i}, a_{i}]$, where $x_{i}\in X$, $a_{i}\in A$.  
Then $[A, A]=[X, A]$.
\end{Lemma}
 
\begin{proof} Let us apply the identity (\cite{Mo}, p.561)
$$
[ab, c]=[a, bc]+[b, ca]. \eqno{(1.1.1)}
$$

Namely, let $a=x_1\dots x_n$; let us induct on $n$ to prove that $[a, 
A]\subset [X, A]$.  For $n=1$ the 
statement is obvious.  If $n>1$, then $a=xa_1$, where $x\in X$ and due 
to (1.1.1) we have
$$
[a, c]=[xa_{1}, c]=[x, a_{1}c]+[a_{1}, cx]
$$
\end{proof}
 
\ssbegin{1.2}{Lemma} Let $A$ be an associative algebra and $a\mapsto 
\omega(a)$ be its involutive antiautomorphism (transposition for $A=\Mat(n)$).  
Let $L$ be an invariant functional on $A$ (like trace, i.e., $L([A, 
A])=0$) such that $L(\omega(a))=L(a)$ for any $a\in A$.  Define the 
bilinear form on $A$ by setting
$$
\langle u, v\rangle =L(u\omega(v)) \text{ for any } u, v\in A. \eqno{(1.2.1)}
$$
Then 

{\em i)} $\langle u, v\rangle =\langle v, u\rangle $;

{\em ii)} $\langle xu, v\rangle =\langle u, \omega(x)v\rangle $;

{\em iii)} $\langle ux, v\rangle =\langle u, v\omega(x)\rangle $;

{\em iv)} $\langle [x, u], v\rangle =\langle u, [\omega(x), 
v]\rangle$.
\end{Lemma}
 
\begin{proof} (Clearly, iii) is similar to ii)).
$$
\renewcommand{\arraystretch}{1.4}
\begin{array}{ll}
\text{ i)}& \langle u, v\rangle =L(u\omega(v)) =L(\omega(u\omega(v)))= 
L(v\omega(u))=\langle v, u\rangle .\\
\text{  ii)} &\langle xu, v\rangle =L(xu\omega(v)) =L(u\omega(v)x)= 
L(u\omega(\omega(x)v))=\langle u, \omega(x)v\rangle .\\
\text{ iv)} &\langle [x, u], v\rangle =\langle xu, v\rangle -
\langle ux, v\rangle\\
&\langle u, \omega(x)v\rangle -\langle u, v\omega(x)\rangle=\langle 
u, [\omega(x), v]\rangle.
\end{array}
$$
\end{proof}

\ssec{1.3.  Traces and forms on $U(\fg)$} Let $\fg$ be a finite 
dimensional Lie algebra, $Z(\fg)$ the center of $(U(\fg)$, $W$ the 
Weyl group of $\fg$ and $\fh$ a Cartan subalgebra of $\fg$.  The 
following statements are proved in \cite{Di}.

\ssbegin{1.3.1}{Proposition} {\em i)} $U(\fg)=Z(\fg)\oplus [U(\fg), U(\fg)]$.
    
{\em ii)} Let $\sharp: Z(\fg)\oplus [U(\fg), U(\fg)]\tto Z(\fg)$ be 
the natural projection.  Then
$$
(uv)^\sharp=(vu)^\sharp \text{ and }(zv)^\sharp=z(v)^\sharp 
\text{ for any } u, v\in U(\fg)\text{ and } z\in Z(\fg).
$$

{\em iii)}  $U(\fg)=S(\fh)^W\oplus [U(\fg), U(\fg)]$.

{\em iv)} Let $\lambda$ be the highest weight of the irreducible 
finite dimensional $\fg$-module $L^\lambda$ and $\varphi$ the 
Harish-Chandra homomorphism.  Then
$$
\varphi(u^\sharp)(\lambda)=\frac{\tr(u|_{L^{\lambda}})}{\dim L^\lambda}.
$$
\end{Proposition}

\ssec{1.3.2} On $U(\fg)$, define a form with values in $Z(\fg)$ by setting
$$
\langle u, v\rangle= (u\omega(v))^\sharp, \eqno{(*)}
$$
where $\omega$ is the Chevalley involution in $U(\fg)$.

\begin{Lemma} The form $(*)$ is nondegenerate on $U(\fg)$. 
\end{Lemma}

\begin{proof} Let $\langle u, v\rangle=0$ for any $v\in U(\fg)$. By 
Proposition 1.3.1
$$
\tr(u\omega(v))=\varphi((u\omega(v))^\sharp)(\lambda)\cdot\dim 
L(\lambda)=\varphi(\langle u, v\rangle)(\lambda)\cdot\dim 
L(\lambda)=0;
$$
hence, $u=0$ on $L(\lambda)$ for any irreducible finite dimensional 
$L(\lambda)$, and, therefore, $u=0$ in $U(\fg)$.
\end{proof}

\ssbegin{1.3.3}{Lemma} For any  $\lambda\in\fh^*$ define a 
$\Cee$-valued form on $U(\fg)$ by setting
$$
\langle u, v\rangle_\lambda=\varphi(\langle u, v\rangle)(\lambda).
$$
The kernel of this form is a maximal ideal in $U(\fg)$. \end{Lemma}

\begin{proof} The form $\langle \cdot, \cdot\rangle_\lambda$ arises 
from a linear functional $L(u)=\varphi(u^\sharp)(\lambda)$; hence, 
by Lemma 1.2 its kernel is a twosided ideal $I$ in $U(\fg)$. On 
$A=U(\fg)/I$, the form induced is nondegenerate. If $z\in Z(\fg)$, 
then
$$
\langle z, v\rangle_\lambda=L(z\omega(v))=L(z)L(\omega(v));
$$
hence, $z-L(z)\in I$. Therefore, the only $\fg$-invariant elements in 
$A$ are those from $\Span(1)$.

Let $J$ be a twosided nontrivial ($\neq A, 0$) ideal in $A$ and 
$J=\mathop{\oplus}\limits_{\mu}J^\mu$ be the decomposition into 
irreducible finite dimensional $\fg$-modules (with respect to the 
adjoint representation). Since $J\neq A$, it follows that $J^0=0$. 
Hence, $L(J)=0$ and $\langle J, A\rangle$. Thus, $J=0$.
\end{proof}

\ssec{1.4. Gelfand--Tsetlin basis and transvector algebras} (For 
recapitulation on transvector algebras see \cite{Zh}.) 

Let $E_{ij}$ be the matrix units.  In $\fgl(3)$, we fix the subalgebra 
$\fgl(2)$ embedded into the left upper corner and let $\fh$ denote the 
Cartan subalgebra of $\fgl(3)=\Span(E_{ii}: i= 1, 2, 3)$.

There is a one-to-one correspondence between finite dimensional 
irreducible representations of $\fgl(3)$ and the sets
$$
(\lambda_{1}, \lambda_{2}, \lambda_{3})\text{ such that } 
\lambda_{1}-\lambda_{2}, \lambda_{2}-\lambda_{3}\in\Zee_{+}.
$$
Such sets are called highest weights of the corresponding irreducible 
representation whose space is denoted $L^\lambda$. With each such 
$\lambda$ we associate a Gelfand--Tsetlin diagram $\Lambda$:
$$
\begin{matrix}
\lambda_{31}&&\lambda_{32}&&\lambda_{33}\cr   
&\lambda_{21}&&\lambda_{22}&&\cr   
&&\lambda_{11}&&\cr   
\end{matrix},\eqno{(1.4.1)}
$$
where the upper line coincides with $\lambda$ and where 
``betweenness'' conditions hold:
$$
\lambda_{k, i}-\lambda_{k-1, i}\in\Zee_{+};\quad  
\lambda_{k-1, i}-\lambda_{k, i+1}\in\Zee_{+}\text{ for any } i=1, 2;\; 
k=2, 3.    \eqno{(1.4.2)}
$$
Set
$$
\begin{matrix}
    z_{21}=E_{21}, \; z_{12}=E_{12};\;  z_{13}=E_{13}, \; z_{32}=E_{32};\cr
z_{31}=(E_{11}-E_{22}+2)E_{31}, \quad z_{23}=(E_{11}-E_{22}+2)E_{23}-
E_{21}E_{13}.
\end{matrix}    \eqno{(1.4.3)}
$$
Set $(L^\lambda)^+=\Span(u: u\in L^\lambda, E_{12}u=0)$.

\ssbegin{1.4.1}{Theorem} {\em (see \cite{Mo})} Let $v$ be a nonzero 
highest weight vector in $L^\lambda$, and $\Lambda$ a Gelfand--Tsetlin 
diagram. Set
$$
v_{\Lambda}=z_{21}^{\lambda_{21}-\lambda_{11}} 
z_{31}^{\lambda_{31}-\lambda_{21}}z_{32}^{\lambda_{32}-
\lambda_{22}}v
$$
and let $l_{ki}=\lambda_{ki}-i+1$. Then

{\em i)} The vectors $v_{\Lambda}$ parametrized by Gelfand--Tsetlin 
diagrams form a basis in $L^\lambda$.

{\em ii)} The $\fgl(3)$-action on vectors $v_{\Lambda}$ is given by the 
following formulas
$$
\renewcommand{\arraystretch}{1.4}
\begin{array}{l}
    E_{11}v_{\Lambda}=\lambda_{11}v_{\Lambda};\\
    E_{22}v_{\Lambda}=(\lambda_{21}+\lambda_{22}-\lambda_{11})v_{\Lambda};\\
    E_{33}v_{\Lambda}=(\mathop{\sum}\limits_{i=1}^3\lambda_{3i}-
    \mathop{\sum}\limits_{j=1}^2\lambda_{2j})v_{\Lambda};\\
    E_{12}v_{\Lambda}=-(l_{11}-l_{21})(l_{11}-l_{22})
    v_{\Lambda+\delta_{11}};\\
    E_{21}v_{\Lambda}=v_{\Lambda-\delta_{11}};\\
 E_{23}v_{\Lambda}=-\displaystyle\frac{(l_{21}-l_{31})(l_{21}-l_{32})(l_{21}-l_{33})}
 {(l_{21}-l_{22})}v_{\Lambda+\delta_{11}}-
 \displaystyle\frac{(l_{22}-l_{31})(l_{22}-l_{32})(l_{22}-l_{33})}
 {(l_{22}-l_{21})}v_{\Lambda+\delta_{22}};\\
 E_{32}v_{\Lambda}=\displaystyle\frac{(l_{21}-l_{11})}
 {(l_{21}-l_{22})}v_{\Lambda-\delta_{11}}+
\displaystyle \frac{(l_{22}-l_{11})}
 {(l_{22}-l_{21})}v_{\Lambda-\delta_{22}},
\end{array}
$$
where $\Lambda\pm\delta_{ki}$ is obtained from $\Lambda$ by 
replacing $\lambda_{ki}$ with $\lambda_{ki}\pm 1$ and we assume that 
$v_{\Lambda}=0$ if $\Lambda$ does not satisfy conditions on 
GTs-diagrams.

{\em iii)} The vectors $v_{\Lambda}$ corresponding to the 
GTs-diagrams with $\lambda_{21}=\lambda_{11}$ form a basis of 
$(L^\lambda)^+$.
\end{Theorem}

\ssec{1.5.1.  Twisted generalized Weyl algebras} Recall definition of 
twisted generalized Weyl algebras introduced in \cite{MT}.  Let $R$ be 
a commutative algebra over $\Cee$, $\Gamma$ a finite nonoriented tree 
with $\Gamma_{0}$ being the set of its vertices and $\Gamma_{1}$ that 
of edges.  Let also $\{\sigma_{i}: i\in\Gamma_{0}\}$ be a set of 
pairwise commuting automorphisms of $R$ and $\{t_{i}: 
i\in\Gamma_{0}\}$ the set of nonzero elements from $R$ satisfying the 
following conditions:
$$
\renewcommand{\arraystretch}{1.4}
\begin{array}{l}
    t_{i}t_{j}=\sigma_{i}^{-1}(t_{j})\sigma_{j}^{-1}(t_{i})\text{ if 
    }(i, j)\in\Gamma_{1}, \\
\sigma_{i}(t_{j})=t_{j}\text{ if }(i, j)\not\in\Gamma_{1}.
\end{array}\eqno{(1.5.1)}
$$

\begin{Lemma} {\em (\cite{MT})} Let $\fA'$ be the algebra generated 
by $X_{i}, Y_{i}: i\in\Gamma_{0}$ subject to the relations (for any $r\in R$)

{\em 1)} $X_{i}r=\sigma_{i}(r)X_{i}$,

{\em 2)} $Y_{i}r=\sigma_{i}^{-1}(r)Y_{i}$,

{\em 3)} $X_{i}Y_j=Y_jX_{i}$ if $i\neq j$,

{\em 4)} $Y_{i}X_{i}=t_{i}$,

{\em 5)} $X_{i}Y_{i}=\sigma_{i}(t_{i})$.

Then $\fA'\neq 0$, $\fA'$ is $\Zee$-graded, and among homogeneous 
(with respect to the grading) twosided ideals of $\fA'$ whose 
intersection with $R$ is trivial is a maximal one, $I$.  \end{Lemma}

The quotient $\fA=\fA'/I$ 
is called the {\it twisted generalized Weyl algebra}.

\ssec{1.5.2. Example} Let $\gamma$ be the Dynkin graph for the root 
system $A_{n-1}$. Let $V$ be an $n$-dimensional vector space with 
basis $e_{1}, \dots , e_{n}$, and let $\eps_{1}, \dots , \eps_{n}$ be 
the dual basis of $V^*$. Let $T=\{\eps_{i}-\eps_{j}\}$ be the root 
system of $A_{n-1}$ in $V^*$ and $\alpha_{1}, \dots , \alpha_{n-1}$ the 
system of simple roots. Set $R=S(V)$ and for $h\in V$ define:
$$
\sigma_{i}(h)=h-\alpha_{i}(h) \text{ for } i= 1, 2, \dots, n-1
$$
and having extending $\sigma_{i}$ to an automorphism of $R$.  Clearly, 
the $\sigma_{i}$ pairwise commute for $i\in \Gamma_{0}$.  In $V$, 
select vectors $h_{1}, \dots , h_{n-1}$ such that 
$\alpha_{i}(h_{j})=\delta_{ij}$ and for an arbitrary collection 
$f_{1}, \dots , f_{n-1}$ of polynomials in one indeterminate set
$$
\renewcommand{\arraystretch}{1.4}
\begin{array}{l}
    t_{1}(v)=f_{1}(h_{1})f_{2}(h_{2}-h_{1}),\quad 
t_{2}(v)=f_{2}(h_{2}-h_{1}+1)f_{3}(h_{3}-h_{2}),\dots , \\ 
t_{n-1}(v)=f_{n-1}(h_{n-1}-h_{n-2}+1)f_{n}(h_{n-1}).
\end{array}$$

It is not difficult to verify that conditions (1.5.1) are satisfied.

\section*{\S 2. Formulations of main results}

\ssec{2.1. Modules $S^\alpha(V)$} Let $\fg=\fgl(3)$ be the Lie 
algebra of $3\times 3$ matrices over $\Cee$. For any $\alpha\in \Cee$ 
denote by $S^\alpha(V)$ the irreducible $\fg$-module with highest 
weight $(\alpha, 0, 0)$.

If $\alpha\in \Zee_{+}$, then $S^\alpha(V)$ is the usual $\alpha$-th 
symmetric power of the identity $\fg$-module $V$. Namely:
$$
S^\alpha(V)=\Span(x_{1}^{k_{1}}x_{2}^{k_{2}}x_{3}^{k_{3}}: 
k_{1}+k_{2}+k_{3}=\alpha; \; \; k_{1}, k_{2}, k_{3}\in \Zee_{+}).
$$
For $\alpha\not\in \Zee_{+}$ we have (like in semi-infinite cohomology 
of Lie superalgebras)
$$
S^\alpha(V)=\Span(x_{1}^{k_{1}}x_{2}^{k_{2}}x_{3}^{k_{3}}: 
k_{1}+k_{2}+k_{3}=\alpha; \; \; k_{2}, k_{3}\in \Zee_{+}).
$$

\begin{Remark} The expression $x^{k}$ for $k\in \Cee$ is 
understood as a formal one, satisfying $\pderf{x^{k}}{x}=kx^{k-1}$. \end{Remark}

On $S^\alpha(V)$ the $\fg=\fgl(3)$-action is given by $E_{ij}\mapsto 
x_{i}\pder{x_{j}}$.

\ssbegin{2.2}{Theorem} {\em i)} $S^\alpha(V)$ is an irreducible 
$\fg$-module for any $\alpha$.

{\em ii)} The kernel $J^\alpha$ of the corresponding to $S^\alpha(V)$ 
representation of 
$U(\fg)$ is a maximal ideal if $\alpha\not\in\Zee_{<0}$.

Set $\fA^\alpha=U(\fg)/J^\alpha$ and let $\theta$ be the highest 
weight of the adjoint representation of $\fg$. Now consider 
$\fA^\alpha$ as $\fg$-module with respect to the adjoint 
representaiton.

{\em iii)} $\fA^\alpha=\mathop{\oplus}\limits_{k=0}^\infty 
L^{k\theta}$ if $\alpha\not\in\Zee_{\geq 0}$.

{\em iv)} $\fA^\alpha=\mathop{\oplus}\limits_{k=0}^\alpha L^{k\theta}$ 
if $\alpha\in\Zee_{\geq 0}$.  

{\em v)} The algebra $\fA^\alpha$ for $\alpha \not\in\Zee_{\geq 0}$ 
is isomorphic to the following twisted generalized Weyl algebra
$$
\renewcommand{\arraystretch}{1.4}
\begin{array}{l}
    R=\Cee[E_{11}, E_{22}, E_{33}]/(E_{11}+ E_{22}+E_{33}-\alpha),\\
    \sigma_1(E_{11})= E_{11}-1,\quad \sigma_1(E_{22})=E_{22}+1,\quad 
    \sigma_1(E_{33})=E_{33},\\
   \sigma_2(E_{11})= E_{11},\quad \sigma_2(E_{22})=E_{22}-1,\quad 
    \sigma_2(E_{33})=E_{33}+1,\\
\end{array}
$$

{\em vi)} The form $\langle u, 
v\rangle_\alpha=\varphi(u\omega(v)^\sharp)(\alpha, 0, 0)$ is 
nondegenerate on $\fA^\alpha$ for $\alpha \not\in\Zee_{< 0}$.
\end{Theorem}

\ssec{2.3} Let $\fh=\Span(E_{11}, E_{22}, E_{33})$ be Cartan 
subalgebra in $\fg$ and $\eps_{1}, \eps_{2}, \eps_{3}$ the dual basis 
of $\fh^*$.  Let $Q=\{\sum k_{i}\eps_{i} : \sum k_{i}=0\}$ be the root 
lattice of $\fg$.  For any $\mu\in Q$ define
$$
(\fA^\alpha)_{\mu}=\{u\in \fA^\alpha: [h, u]=\gamma(h)u\text{ for any 
}h\in\fh\}. \eqno{(2.3.1)}
$$
Clearly, $\fA^\alpha$ is $Q$-graded:
$$
\fA^\alpha=\mathop{\oplus}\limits_{\mu\in Q}(\fA^\alpha)_{\mu}.
$$
Theorem 2.4 below shows that $(\fA^\alpha)_{\mu}=Ru_{\mu}$, where
$u_{\mu}\in \fA^\alpha$ is defined uniquely up to a constant factor 
and $R=\Cee[E_{11}, E_{22}, E_{33}]/(E_{11}+E_{22}+ E_{33}-\alpha)$.

Denote by $(\fA^\alpha)^+$ the subalgebra of $\fg$ consisting of 
vectors highest with respect to the fixed $\fgl(2)$:
$$
(\fA^\alpha)^+=\{u\in \fA^\alpha:[E_{12}, u]=0\}. \eqno{(2.3.2)}
$$

The algebra $(\fA^\alpha)^+$ also admits $Q$-grading:
$$
(\fA^\alpha)^+=\mathop{\oplus}\limits_{\nu\in Q}(\fA^\alpha)^+_{\nu}.
\eqno{(2.3.3)}
$$
Denote: $Q^+=\{\nu\in Q: (\fA^\alpha)^+_{\nu}\neq 0$. 

Theorem 2.4 below shows that $(\fA^\alpha)^+_{\nu}=
\Cee[E_{33}]u_{\nu}^+$, where $\nu\in Q^+$. For $f, g\in
\Cee[E_{33}]$ and $\nu\in Q^+$ set
$$
\langle f, g\rangle_{\nu}^+=\langle fu_{\nu}, gu_{\nu}\rangle_{\alpha}. 
\eqno{(2.3.4)}
$$

For $f, g\in R$ and $\mu\in Q$ set
$$
\langle f, g\rangle_{\mu}=\langle fu_{\mu}, gu_{\mu}\rangle_{\alpha}. 
\eqno{(2.3.5)}
$$

For  $k\geq 0$ and $\nu\in Q^+$ set
$$
f_{k, \nu}(E_{33})u_{\nu}=
\begin{cases}(\ad z_{31} )^k(u_{\nu+k(\eps_{1}-\eps_{3})})&\text{ for 
}\nu(E_{33})\leq 0\cr
(\ad z_{23} )^k(u_{\nu+k(\eps_{3}-\eps_{2})})&\text{ for 
}\nu(E_{33})\geq 0\end{cases}\eqno{(2.3.6)\choose (2.3.7)}
$$

For  $k, l\geq 0$ and $\nu\in Q$ set
$$
f_{l, k}^\nu(E_{11}, E_{22}, E_{33})u_{\nu}=
\begin{cases}(\ad z_{21} )^l(\ad z_{31} 
)^k(u_{\nu+k(\eps_{1}-\eps_{3})+l(\eps_{1}-\eps_{2})})&\text{ for 
}\nu(E_{33})\leq 0\cr (\ad z_{21} )^l(\ad z_{23} 
)^k(u_{\nu+k(\eps_{3}-\eps_{2})+l(\eps_{1}-\eps_{2}})&\text{ for 
}\nu(E_{33})\geq 0\end{cases}\eqno{(2.3.8)\choose (2.3.9)}
$$

\ssbegin{2.4}{Theorem} {\em 0)}
$(\fA^\alpha)_{\nu}^+=\Cee[E_{33}]u_{\nu}^+$, where $u_{\nu}$ is
determined uniquely up to a constant factor.

{\em 1)} $\langle (\fA^\alpha)^+_{\nu}, 
\fA^\alpha)^+_{\nu}\rangle_{\alpha}=0$ for $\nu\neq \mu$.

{\em 2)}  The polynomials $f_{k, \nu}(E_{33})$ are orthogonal relative 
$\langle \cdot, \cdot\rangle_{\nu}^+$.

{\em 3)}   The polynomials $f_{k, \nu}(E_{33})$ satisfy the difference 
equation
$$
(E_{33}-\nu(E_{33})+1)(E_{33}+\nu(E_{11})-\alpha)\Delta f-
E_{33}(E_{33}+\nu(E_{22})-\alpha-2)\nabla f=k(k+2\nu(E_{11}+2)f 
\text{ if }\nu(E_{33})<0
$$
$$
(E_{33}+1)(E_{33}+\nu(E_{11})-\alpha)\Delta f- 
(E_{33}-\nu(E_{33}))(E_{33}+\nu(E_{22})-\alpha-2)\nabla 
f=k(k-2\nu(E_{11}+2)f \text{ if }\nu(E_{33})\geq 0
$$

{\em 4)}   Explilcitely, $f_{k, \nu}(E_{33})$ is of the form
$$
f_{k, \nu}(E_{33})=\text{const}\times\, {}_{3}F_{2}\left(\begin{matrix}-k, k+2k_1+2, 
-E_{33}\cr 1-k_3 
k_1-\alpha\end{matrix}~\vert~ 1\right),
$$
where 
$$
{}_{3}F_{2}\left(\begin{matrix}\alpha_{1}, \alpha_{2}, \alpha_{3}\cr 
\beta_{1}, \beta_{2}\end{matrix}\left \vert \right.
z\right)=\mathop{\sum}\limits_{i=0}^\infty
\displaystyle\frac{(\alpha_{1})_{i}(\alpha_{2})_{i} 
(\alpha_{3})_{i}}{(\beta_{1})_{i}(\beta_{2})_{i}}\,
\displaystyle\frac{z^i}{i!}
$$
is a generalized hypergeometric function, $(\alpha)_{0}=1$  and 
$(\alpha)_{i}=\alpha(\alpha+1)\dots(\alpha+i-1)$ for $i>0$. \end{Theorem} 

\ssbegin{2.5}{Theorem} {\em 0)} $(\fA^\alpha)_{\nu}=\Cee[E_{11}, 
E_{22}, E_{33}]u_{\nu}$, 
where $u_{\nu}$ is determined uniquely up to a constant factor.

{\em 1)} $\langle (\fA^\alpha)_{\nu}, 
\fA^\alpha)_{\nu}\rangle_{\alpha}=0$ for $\nu\neq \mu$.

{\em 2)} The polynomials $f_{l, k}^\nu(E_{11}, E_{22}, E_{33})$ form 
an orthogonal basis of $R$ relative $\langle \cdot, 
\cdot\rangle_{\nu}$.

{\em 3)} The polynomials $w(f_{l, k})(E_{11}, E_{22}, E_{33})$ for 
$w\in W$ form an orthogonal basis of $R$ relative $\langle \cdot, 
\cdot\rangle_{w(\nu)}$ provided polynomials $f_{l, k}(E_{11}, E_{22}, 
E_{33})$ form an orthogonal basis of $R$ relative $\langle \cdot, 
\cdot\rangle_{\nu}$.

{\em 4)} The polynomials $f_{l, k}^\nu(E_{11}, E_{22}, E_{33})$ for 
$\nu\in Q^+$ and $\nu(E_{33})\leq 0$ satisfy the system of two 
difference equations
$$
\renewcommand{\arraystretch}{1.4}
\begin{array}{l}
    {}[f(H_{1}+2, H_{2})-f(H_1, H_{2})]\cdot 
    \frac14(H_{1}-H_{2}+\alpha+1)(H_{1}+H_{2}-\alpha)-\\
{}[f(H_{1}, H_{2})-f(H_1-2, H_{2})]\cdot 
\frac14(H_{1}-H_{2}+\alpha-\nu(E_{11}))(H_{1}+H_{2}-
\alpha-1+\nu(E_{22}))=\\
{}[l^2+l(\nu(E_{11})+\nu(E_{22})+1))+\nu(E_{22})-\nu(E_{11})]f,\\
{}\\
{}[2\alpha-\nu(H_{2})(\alpha+2+\nu(H_{2}))+H_{2}(2\alpha+1+2\nu(H_{2}))-
2H_{2}^2]f(H_{1}, H_{2})-\\
\frac12(H_{2}+1-\nu(H_{2}))(H_{1}-H_{2}+\alpha-2\nu(E_{11}))
f(H_{1}-1, H_{2}+1)-\\
\frac12H_{2}(H_{1}-H_{2}+\alpha+2)f(H_{1}+1, H_{2}+1)-\\
\frac12(H_{2}+1-\nu(H_{2}))(\alpha-H_{1}-H_{2})f(H_{1}+1, H_{2}+1)-\\
\frac12H_{2}(\alpha-H_{1}-H_{2}+2-2\nu(E_{22}))f(H_{1}-1, H_{2}-1)=\\
{}[2k^2+4kl+4k(1+\nu(E_{11}))+2l(1+\nu(E_{11})-\nu(E_{22}))+\\
\nu(E_{11})^2-\nu(E_{22})^2+4\nu(E_{11})]f(H_{1}, H_{2}).\end{array}
$$
\end{Theorem}

\section*{\S 3. Proof of Theorem 2.2}

i) The module $S^\alpha(V)$ is irreducible if and only if it has no 
vacum vectors (i.e, vectors annihilated by $E_{12}$ and $E_{23}$.  
This is subject to a direct verification.

ii) Follows from Excercise 858 of Ch. 8 of \cite{Di}.

iii) Let $A_{3}$ be the Weyl algebra (i.e., it is generated by the 
$p_{i}$ and $q_{i}$ for $i= 1, 2, 3$ satisfying 
$$
p_{i}p_{j}-p_{j}p_{i}=q_{i}q_{j}-q_{j}q_{i}=0;\quad 
p_{i}q_{j}-q_{j}p_{i}=-\delta_{ij}.\eqno{(3.1)}
$$
Setting $E_{ij}\mapsto p_{i}q_{j}$ we see that the homomorphism 
$\varphi: U(\fg)\tto\End(S^\alpha(V))$ factors through $A_{3}$ and 
$A_{3}$ acts on $S^\alpha(V)$ so that $p_{i}\mapsto x_{i}$ and $q_{i}\mapsto 
\pder{x_{i}}$. Let us describe the image of $\varphi$. To this end, 
on $A_{3}$, introcude a grading by setting
$$
\deg p_{i}= 1\quad \deg q_{i}= -1\text{ for $i= 1, 2, 3$ }.\eqno{(3.2)}
$$
Now it is clear that $\IM \varphi$ is the algebra $B_{3}$ of elements 
of degree 0. 

To describe highest weight elements in $B_{3}$, it suffices to 
describe same in $S^k(V)\otimes S^k(V^*)$. Let us identify $S^k(V)\otimes 
S^k(V^*)$ with $\End(S^k(V))$, let $u\in\End(S^k(V))$ commutes with 
the action of $E_{12}$ and $E_{23}$ on $S^k(V)$. But then $u$ is 
uniquely determined by its value on the lowest weight vector 
$x_{3}^k\in S^k(V)$; moreover, $E_{12}x_{3}^k=0$. Hence,
$$
u(x_{3}^k)=a_{0}x_{3}^k+\mathop{\sum}\limits_{i=0}^k
a_{i}x_{1}^{i}x_{3}^{k-i},
$$
so
$$
u(x_{3}^k)=\frac{1}{k}a_{0}(\mathop{\sum}\limits_{i=0}^kx_{i}\pder{x_{i}})x_{3}^k+\mathop{\sum}\limits_{i=0}^k
\frac{(k-i)!}{k!}a_{i}(x_{1}\pder{x_{3}})^ix_{3}^{k}.
$$
This shows that the algebra of highest weight vectors in $B_{3}$ is 
generated by $p_{1}q_{3}$ and $z=p_{1}q_{1}+p_{2}q_{2}+p_{3}q_{3}$. 
If $\alpha\not\in \Zee_{\geq 0}$, then $\fA_{\alpha}$ is the quotient 
of $B_{3}$ modulo $(z-\alpha)$. This proves iii).

iv) In this case $\fA_{\alpha}=\End(S^k(V))$ and the proof follows 
from the arguments at the end of the above paragraph.

v) For the canonical homomorphism $\varphi: U(\fg)\tto \fA_{\alpha}$ 
and any $H\in R$ set
$$
\renewcommand{\arraystretch}{1.4}
\begin{array}{l}
    X_{1}=\varphi(E_{12}), \; X_{2}=\varphi(E_{23}), \; 
Y_{1}=\varphi(E_{21}), \; Y_{2}=\varphi(E_{32}); \\
R=\Cee[E_{11}, E_{22}, E_{33}]/(z-\alpha);\\
\sigma_{1}(H)=H-(\eps_{1}-\eps_{2})(H),\; 
\sigma_{2}(H)=H-(\eps_{2}-\eps_{3})(H);\\
t_{1}=E_{22}(E_{11}+1),\; 
t_{2}=E_{33}(E_{22}+1).
\end{array}
$$
Now, as is easy to verify, all the relations of sec. 1.5.1 are 
satisfied. If $I$ is a twosided ideal of $\fA_{\alpha}$, then thanks 
to iii) and iv) it contains $L^{k\omega}$ for some $k\geq 0$ and, 
therefore, elements of weight 0 with respect to the adjoint action of 
Cartan subalgebra of $\fg$. Hence, $I\cap R\neq 0$. Hence, there 
exists a surjection $\psi: \fA'\tto\fA_{\alpha}$.

vi) By 1.3.3 the kernel of $\langle\cdot, \cdot\rangle_{\alpha}$ in 
$U(\fg)$ is a maximal ideal. But $\fA_{\alpha}=U(\fg)/J^\alpha$, where
$J^\alpha$ is maximal due to i). So $J^\alpha$ coincides with the 
kernel of $\langle\cdot, \cdot\rangle_{\alpha}$ in 
$U(\fg)$ and the form is nondegenerate on $\fA_{\alpha}$.

\section*{\S 4. Proof of Theorem 2.4}

0) Direct computations show that the set of elements from $A_3$ 
commuting with $E_{12}$ is a subalgebra generated by $p_{1}$, 
$q_{2}$, $p_{3}$, $q_{3}$ and $z=p_{1}q_{1}+p_{2}q_{2}+p_{3}q_{3}$. So 
this algebra is the linear span of the elements of the form
$$
u=p_{1}^{k_{1}}q_{2}^{k_{2}}p_{3}^{k_{3}}q_{3}^{k_{4}}z^{k_{5}}.
$$
If $u\in B_{3}$, then $k_{1}+k_{3}=k_{2}+k_{4}$, so
$$
u=\begin{cases}p_{1}^{k_{1}}q_{2}^{k_{2}}p_{3}^{k_{3}-k_{4}}
p_{3}^{k_{4}}q_{3}^{k_{4}}z^{k_{5}}&\text{if }k_{3}\geq k_{4}\\
p_{1}^{k_{1}}q_{2}^{k_{2}}p_{3}^{k_{3}}
q_{3}^{k_{3}}q_{3}^{k_{4}-k_{3}}z^{k_{5}}&\text{if }k_{3}\leq k_{4}.
\end{cases}\eqno{(4.1)}
$$
Hence, setting for $\nu=\sum k_{i}\eps_{i}$ such that $\sum 
k_{i}=0$, $k_1\geq 0$ and $k_2\leq 0$
$$
u_{\nu}^+=\begin{cases}p_{1}^{k_{1}}q_{2}^{-k_{2}}p_{3}^{k_{3}}&\text{if 
}k_{3}\geq 0\\
p_{1}^{k_{1}}q_{2}^{-k_{2}}
q_{3}^{-k_{3}}&\text{if }k_{3}\leq 0
\end{cases}\eqno{(4.2)}
$$
we obtain the statement desired.

1) Let $u\in(\fA^+_\alpha)_{\mu}$, $v\in(\fA^+_\alpha)_{\nu}$, 
and $h\in\fh$. Then by heading iv) of Lemma 1.2 we obtain:
$$
\langle [h, u], v\rangle =\mu(h)\langle u, v\rangle=\langle u, [h, v]\rangle
=\nu(h)\langle u, v\rangle.
$$
So $\langle u, v\rangle=0$ if $\mu\neq \nu$.

2) Let $\nu(E_{33})\leq 0$. We have:
$$
\renewcommand{\arraystretch}{1.4}
\begin{array}{l}
f_{k, \nu}u_{\nu}^+=(\ad z_{31})^k(u_{\nu+k(\eps_{1}-\eps_{3})}^+)=
\ad z_{31}(\ad 
z_{31})^{k-1}(u_{\nu+(k-1)(\eps_{1}-\eps_{3})+(\eps_{1}-\eps_{3})}^+)=\\
(\ad z_{31})f_{k-1,
\nu+(\eps_{1}-\eps_{3})}u_{\nu+(\eps_{1}-\eps_{3})}^+.
\end{array}
$$

Direct verification shows that
$$
\renewcommand{\arraystretch}{1.4}
\begin{array}{l}
(\ad
z_{31})(fu_{\nu}^+)=\{(E_{33}-\nu(E_{33}))[(E_{33}-\alpha)\nu(E_{11})+
(\nu(E_{11})-1)\nu(E_{11})]f(h)-\\
E_{33}[(E_{33}-\alpha)\nu(E_{11})-(\nu(E_{22})+2)\nu(E_{11})]f(h-1)\}
u_{\nu-(\eps_{1}-\eps_{3})}.
\end{array}\eqno{(4.3)}
$$
It easily follows from Lemma 1.2 that for any $z\in U(\fg)$ we have
$$
\langle (\ad z)(u), v\rangle=\langle u, (\ad \omega(z))(v)\rangle,
$$
but 
$$
\omega(z_{31})=\omega((E_{11}-E_{22}+2)E_{31}+E_{21}E_{32})=
E_{13}(E_{11}-E_{22}+2)+E_{23}E_{12}.
$$

Since $fu_{\nu}$ is a highest weight vector with respect to the fixed 
$\fgl(2)$, it follows that
$$
\renewcommand{\arraystretch}{1.4}
\begin{array}{l}
    (\ad\omega(z_{31}))(fu_{\nu})=(\ad 
(E_{13}(E_{11}-E_{22}+2))(fu_{\nu})=\\
(\nu(E_{11})-\nu(E_{22})+2)\Delta f\cdot u_{\nu+(\eps_{1}-\eps_{3})}.
\end{array}
$$
Now, let us induct on $k$. For $k=0$ the statement is obvious. For 
$k>0$ and $\deg g<k$ we have
$$
\renewcommand{\arraystretch}{1.4}
\begin{array}{l}
   \langle f_{k, \nu}, g\rangle_{\nu}=\langle f_{k, \nu}u_{\nu}^+, 
   gu_{\nu}^+\rangle=\\
\langle f_{k-1, \nu+(\eps_{1}-\eps_{3})}u_{\nu+(\eps_{1}-\eps_{3})}^+, 
   (\ad(\omega(z_{31})))(gu_{\nu}^+)\rangle=\\
   \langle f_{k-1, \nu+(\eps_{1}-\eps_{3})}, 
  (\nu(E_{11})-\nu(E_{22})+2)\Delta 
  g\rangle_{\nu+(\eps_{1}-\eps_{3})}=0
\end{array}
$$
by inductive hypothesis.

The case $\nu(E_{33})\geq 0$ is similar. \qed

3) Observe that $z=E_{13}E_{31}+E_{23}E_{32}$ belongs to the 
centralizer of $\fgl(2)$ in $U(\fg)$. Let $\nu(E_{33})\leq 0$. Then 
$u_{\nu}^+=p_{1}^{k_{1}}q_{2}^{k_{2}}q_{3}^{k_{3}}$ as in (4.1.2). 
Having applied $\ad z$ to $fu_{\nu}^+$ we obtain:
$$
\renewcommand{\arraystretch}{1.4}
\begin{array}{l}
(\ad z)(fu_{\nu}^+)=E_{13}E_{31}fu_{\nu}^++fu_{\nu}^+E_{13}E_{31}- 
E_{13}fu_{\nu}^+E_{31}-E_{31}fu_{\nu}^+E_{13}+\\
E_{23}E_{32}fu_{\nu}^++fu_{\nu}^+E_{32}E_{23}-E_{23}fu_{\nu}^+E_{32}-
E_{32}fu_{\nu}^+E_{23}=\\
E_{11}(E_{33}+1)fu_{\nu}^++fu_{\nu}^+E_{33}(E_{11}+1)-f(E_{33}+1)
u_{\nu}^+E_{11}(E_{33}+1)-
f(E_{33}-1)E_{33}(E_{11}+1)u_{\nu}^++\\
E_{22}(E_{33}+1)fu_{\nu}^++fu_{\nu}^+E_{33}(E_{22}+1)-f(E_{33}+1)
E_{22}u_{\nu}^+(E_{33}+1)-
f(E_{33}-1)E_{33}u_{\nu}^+(E_{22}+1)=\\
(E_{11}+E_{22})(E_{33}+1)fu_{\nu}^++(E_{33}-\nu(E_{33}))(E_{11}+1-
\nu(E_{11})+E_{22}+1-\nu(E_{22}))fu_{\nu}^+-\\
f(E_{33}+1)\cdot(E_{33}+1-\nu(E_{33}))(E_{11}+E_{22}-\nu(E_{11}))u_{\nu}^+-
f(E_{33}-1)E_{33}(E_{11}+E_{22}-\nu(E_{22})+2)u_{\nu}^+=\\
{}[f(E_{33}+1)\cdot(E_{33}+1-\nu(E_{33}))(E_{33}-\alpha+\nu(E_{11}))
+f(E_{33}-1)E_{33}
(E_{33}+\nu(E_{22})-\alpha-2)-\\
(E_{33}-\alpha)(E_{33}+1)f  - 
(E_{33}-\nu(E_{33}))(E_{33}+\nu(E_{11})+
\nu(E_{22})-\alpha-2)f]u_{\nu}^+.
\end{array}
$$
This gives us the right hand side of the first equation of heading 3).

Since $\ad z$ commutes with the $\fgl(2)$-action and preserves the 
degree of polynomial $f$, it follows that $(\ad z)(fu_{\nu})=c\cdot 
(fu_{\nu})$. Counting the constant factor, we arrive to the first 
equation of heading 3).

The proof of the second equation is similar.

\section*{\S 5. Proof of Theorem 2.5}

0) Recall that $B_3$ is the subalgebra of $A_3$ of the elements of 
degree 0 relative grading (3.2).

For $k\in\Zee$ set $r_{i}^k=\begin{cases}=p_{i}^k&\text{if }k\geq 0\cr, 
q_{i}^{-k}&\text{if }k\leq 0\end{cases}$. For
$\gamma=\sum k_{i}\eps_{i}$, where $\sum k_{i}=0$, set
$$
u_{\gamma}=r_{1}^{k_{1}}r_{2}^{k_{2}}r_{3}^{k_{3}}.
$$
Clearly, $B_3$ is the linear span of the elements of the form
$$
p_{1}^{m_{1}}q_{1}^{l_{1}}p_{2}^{m_{2}}q_{2}^{l_{2}}p_{3}^{m_{3}}q_{3}^{l_{3}}, 
\text{ where }m_{1}+m_{2}+m_{3}=l_{1}+l_{2}+l_{3}.
$$
It is also clear teat each such element can be represented in the form
$$
f(E_{11}, E_{22}, E_{33})r_{1}^{k_{1}}r_{2}^{k_{2}}r_{3}^{k_{3}}.
$$
This completes the proof of heading 0).

1) Proof is similar to that from sec. 4.2.

2) Let $\nu(E_{33})\leq 0$. By setting $H_{1}=E_{11}- E_{22}$, 
$H_{2}=E_{22}- E_{33}$ we identify $R=\Cee[E_{11}, E_{22}, 
E_{33}]/(E_{11}+ E_{22},+E_{33}-\alpha)$ with $\Cee[H_{1}, H_{2}]$. 
Let $\Lambda$ is a Gelfand--Tsetlin diagram of the following form:
$$
\begin{matrix}
\nu(E_{11})+k+l&&0&&-(\nu(E_{11})+k+l)\cr   
&\nu(E_{11})+l&&-(\nu(E_{22})+l)&&\cr   
&&\nu(E_{11})&&\cr   
\end{matrix},\eqno{(5.0)}
$$
From the explicit formula for $f_{k, l}^\nu$ we derive that 
$$
f_{k, l}^\nu u_{\nu}=v_{\Lambda}.\eqno{(5.1)}
$$

Now, consider the following operators from the maximal commutative 
subalgebra of $U(\fg)$:
$$
\renewcommand{\arraystretch}{1.4}
\begin{array}{l}
    E_{11}, E_{22}, \\
    \Omega_{2}=E_{11}^2+E_{22}^2+E_{11}-E_{22}+ 2E_{21}E_{12},\\
    \Omega_{3}=E_{11}^2+E_{22}^2+E_{33}^2+E_{11}-E_{22}+ E_{11}-E_{33}+
    E_{22}-E_{33}+\\
    2E_{21}E_{12}+2E_{31}E_{13}+2E_{32}E_{23}.
\end{array}\eqno{(5.2)}
$$
Then we have:
$$
\renewcommand{\arraystretch}{1.4}
\begin{array}{l}
E_{11}v_{\Lambda}=\nu(E_{11})v_{\Lambda}; \; 
E_{22}v_{\Lambda}=-\nu(E_{22})v_{\Lambda};\\
\Omega_{2}v_{\Lambda}=[2l^2+2l(\nu(E_{11})+\nu(E_{22})+1)+\\
\nu(E_{11})^2+\nu(E_{22})^2]v_{\Lambda};\\
\Omega_{3}v_{\Lambda}=2(\nu(E_{11})+k+l)(\nu(E_{11})+k+l+2]v_{\Lambda}.
\end{array}\eqno{(5.3)}
$$
It is easy to check that the operators (5.2) satisfy
$$
\omega(E_{11})=E_{11};\; \omega(E_{22})=E_{22};\; 
\omega(\Omega_{2})=\Omega_{2};\; \omega(\Omega_{3})=\Omega_{3}
$$
and, therefore, they are selfadjoint relative the form $\langle \cdot, 
\cdot\rangle$.  Formula (5.3) makes it manifest that operators (5.2) 
separate the vecotrs $v_{\Lambda}$, hence, these vectors are pairwise 
orthogonal.  Moreover, it is easy to see that $f_{k, l}^\nu$ is of the 
form
$$
f_{k, l}^\nu=H_{1}^lH_{2}^k+\dots,
$$
where the dots designate the summands of degrees $\leq k+l$ of the 
form $H_{1}^aH_{2}^b$, where $(a, b)<(l, k)$ with respect to the 
lexicographic ordering.  Thus, the $f_{l, k}^\nu$ constitute a basis 
of $\Cee[H_{1}, H_{2}]$.

3) The statement follows from the fact that the Weyl group acts on 
$\fA_{\alpha}$ and preserves the form $\langle \cdot, \cdot\rangle$.

4) Since the polynomials $f_{l, k}^\nu u_\nu$ are elements of a 
Gelfand--Tsetlin basis, they are eigenvectors for 
$\Omega_{2}$ and $\Omega_{3}$ with respect to the 
adjoint action of $\fg=\fgl(3)$ on $\fA_{\alpha}$. As we have shown in 
sec 5.2, we have
$$
\renewcommand{\arraystretch}{1.4}
\begin{array}{l}
    \Omega_{2}f_{l, k}^\nu u_\nu=[2l^2+2l(\nu(E_{11})+\nu(E_{22})+1)+\\
\nu(E_{11})^2+\nu(E_{22})^2]f_{l, k}^\nu u_\nu;\\
\Omega_{3}f_{l, k}^\nu u_\nu=2(\nu(E_{11})+k+l)(\nu(E_{11})+k+l+2]
f_{l, k}^\nu u_\nu.
\end{array}
$$
To derive the corresponding equations, we have to explicitely compute 
the actions of $\Omega_{2}$ and $\Omega_{3}$ on $fu_\nu$. Let 
$(\nu(E_{33})\leq 0$; then
$$
\renewcommand{\arraystretch}{1.4}
\begin{array}{l}
    \Omega_{2}fu_\nu=[\nu(E_{11})^2+\nu(E_{22})^2+\nu(E_{11})-
    \nu(E_{22})]fu_\nu+\\
{}2[f(H_{1}, H_{2})-f(H_{1}+2, H_{2})] \cdot \frac14(H_{1}-H_{2}+\alpha+1)
(\alpha-H_{1}-H_{2})u_\nu+\\
{}2[f(H_{1}, H_{2})-f(H_{1}-2, H_{2})] \cdot \frac14(H_{1}-H_{2}+
\alpha-\nu(E_{11}))
(\alpha-H_{1}-H_{2}+1-\nu(E_{22}))u_\nu;\\
{}\\
(\Omega_{3}-\Omega_{2})fu_\nu=(E_{33}^2-E_{11}-E_{22}+2E_{33}+2 
E_{31}E_{13}+2 E_{32}E_{23})fu_\nu=\\
(E_{33}^2-E_{11}-E_{22}+2E_{33}+E_{13}E_{31}+E_{23}E_{32})fu_\nu=\\
(\nu(E_{33})^2-\nu(E_{11})-\nu(E_{22})+2\nu(E_{33})fu_\nu+\\
(E_{11}+E_{22})(E_{33}+1)f(H_{1}, H_{2})u_\nu+\\
(E_{33}-\nu(E_{33})(E_{11}+E_{22}+2-\nu(E_{11})-\nu(E_{22}))f(H_{1}, 
H_{2})u_\nu-\\
f(E_{11}-1, E_{33}+1)(E_{11}-\nu(E_{11}))(E_{33}+1-\nu(E_{33}))u_\nu-\\
f(E_{11}+1, E_{33}-1)E_{33}(E_{11}+1)u_\nu-\\
f(E_{22}-1, E_{33}+1)E_{22}(E_{33}+1-\nu(E_{33}))u_\nu-\\
f(E_{22}+1, E_{33}-1)E_{33}(E_{22}+1-\nu(E_{22}))u_\nu=\\
{}[(\alpha-H_{2})(H_{2}+1)f(H_{1}, H_{2})+  (H_{2}-\nu(H_{2}))
(\alpha-H_{2}+2+\nu(H_{2}))f(H_{1}, H_{2})-\\
f(H_{1}-1, H_{2}+1)(H_{2}+1-\nu(H_{2}))\frac12
(H_{1}-H_{2}+\alpha-2\nu(E_{11}))-\\
f(H_{1}+1, H_{2}+1)H_{2}\frac12
(H_{1}-H_{2}+\alpha+2)-\\
f(H_{1}+1, H_{2}+1)\frac12(\alpha-H_{1}-H_{2})
(H_{2}+1-\alpha-\nu(H_{2}))-\\
f(H_{1}-1,
H_{2}-1)H_{2}\frac12(\alpha-H_{1}-H_{2}+2-2\nu(E_{22}))]u_\nu.
\end{array}
$$
This implies the second equation.

\end{document}